\documentclass[lettersize,journal]{IEEEtran}
\usepackage{amsmath,amsfonts}
\usepackage{algorithmic}
\usepackage{algorithm}
\usepackage{array}
\usepackage{array,multirow,graphicx}
\usepackage[caption=false,font=normalsize,labelfont=sf,textfont=sf]{subfig}
\usepackage{textcomp}
\usepackage{stfloats}
\usepackage{url}
\usepackage[dvipsnames]{xcolor}
\usepackage{verbatim}
\usepackage{graphicx}
\usepackage{cite}
\usepackage{caption}
\usepackage{nomencl}
\usepackage{amsmath}
\makenomenclature
\begin{document}

\title{An Advanced Fuel Efficiency Optimization Model with Fractional Programming}

\author{Md Isfakul Anam,~\IEEEmembership{~Clarkson University}, Tuyen Vu,~\IEEEmembership{~Clarkson University}
\thanks{Md Isfakul Anam, T. T. Nguyen, and T. Vu with Clarkson University, Potsdam, NY, USA; Emails: isfakum@clarson.edu, tnguyen@clarkson.edu, tvu@clarkson.edu. Corresponding Author: T. Vu, Email: tvu@clarkson.edu}}

\maketitle

\begin{abstract}
Reducing the fuel consumption within a power network is crucial to enhance the overall system efficiency and minimize operating costs. Fuel consumption minimization can be achieved through different optimization techniques where the output power of the generators is regulated based on their individual efficiency characteristics. Existing studies primarily focus either on maximizing the efficiency function or minimizing the operating cost function of the generators to minimize fuel consumption. However, for practical implementation, it becomes imperative to incorporate a function within the optimization framework to represent the fuel consumption rate directly. This study introduces a novel approach by formulating a minimization problem with a sum-of-ratios objective function representing the fuel consumption rate. However, optimization problems with sum-of-ratios objective functions or constraints are extremely challenging to solve because of their strong nonlinearity. To efficiently solve the formulated problem, a fractional programming (FP) approach is adopted in this study. This reformulation technique significantly reduces the solution time of the optimization problem and provides a better solution than nonlinear programming (NLP). In addition, the reformulated problem can also be applied to large-scale systems where the NLP fails to converge. The proposed methodology of this study is tested on the notional MVAC ship system, modified IEEE 30-bus and IEEE 118-bus systems. The results demonstrate that the model successfully minimizes fuel consumption by effectively scheduling the generator and ESS dispatch. 

\end{abstract}


\nomenclature[01]{\(P_{i,g}, Q_{i,g}\)}{Real and reactive power output of $i$-th generator}
\nomenclature[02]{\(P_{i,l}, Q_{i,l}\)}{Real and reactive load at $i$-th bus}
\nomenclature[03]{\(P_{i,inj}, Q_{i,inj}\)}{Real and reactive power injection at $i$-th bus}
\nomenclature[04]{\(V_i, \theta_{ik}\)}{voltage and voltage angle difference}
\nomenclature[05]{\(G_{ik}, B_{ik}\)}{Conductance and susceptance of line $ik$}
\nomenclature[06]{\(P_{ik}\)}{Real power capacity of line $ik$}
\nomenclature[07]{\(DR_i, UR_i\)}{Down and up ramp rate of $i$-th generator}
\nomenclature[08]{\(P_i^{C,max}, P_i^{D,max}\)}{Maximum charging and discharging rate of $i$-th ESS}
\nomenclature[09]{\(E_{i,b}\)}{Energy stored at $i$-th ESS}
\nomenclature[10]{\(P_{i,b}^r\)}{Output power of $i$-th ESS}
\nomenclature[11]{\(SOC_i\)}{State of charge of $i$-th ESS}
\nomenclature[12]{\(\alpha\)}{fuel energy density}
\nomenclature[13]{\(N_g\)}{Number of generator}
\nomenclature[14]{\(T\)}{Total planning horizon}
\nomenclature[15]{\(\triangle t\)}{Time step}
\nomenclature[16]{\(B\)}{Number of buses}

\printnomenclature

\section{Introduction}

Continuous increments in load demand on power systems due to the growing number of consumers impose a great challenge for the engineers and operators in the field. It is essential to supply the increased load requirements of a power system territory by introducing additional energy sources and/or expanding the capacity of the existing sources \cite{mistry_enhancement_2014}. Both solutions result in higher fuel costs incurred for operating a large number of distributed generators. Furthermore, the survival period of stand-alone power networks, such as islanded microgrids or ship electrical power systems, largely depends on the fuel consumption of the distributed generators. As a result, a state-of-the-art system efficiency model is required to optimize the fuel consumption of a system.

System efficiency optimization techniques refer to minimizing the fuel consumption rate by regulating the output power of the generators within the system. Over the past decades, a significant number of research has been conducted to improve system efficiency using various optimization methods. 
One notable contribution is presented in \cite{yao_efficiency_2018}, where the authors propose a novel approach to analyze the characteristics of the efficiency function, leading to the determination of the maximum total power supply and overall efficiency for such systems. In \cite{orero_genetic_1998}, a genetic algorithm modeling framework is presented where the optimization problem involves minimizing the thermal cost function. The cost function is constructed based on the power generation characteristics of the hydropower plant. A similar approach can be found in \cite{yuan_enhanced_2008}, where the daily optimal generation scheduling problem (DOHGSB) is solved by implementing a unique differential evolution algorithm. The objective function is formulated by analyzing the hydropower plant characteristics or input-output curve. 
In \cite{vu_large-scale_2018}, the authors introduce a novel distributed algorithm to maximize system efficiency. A fourth-order and a third-order efficiency function for the main and auxiliary power generation module (PGM), respectively, are optimized using a distributed crow search algorithm (DCSA). 
The approaches in \cite{yao_efficiency_2018}-\cite{vu_large-scale_2018} focus on optimizing the efficiency function or the generation characteristics of the generators, which can be complex to apply to systems consisting of generators with different ratings. Also, since the generator efficiency function is not a straightforward representation of the system's fuel consumption, it is impractical to optimize the efficiency function with an objective to improve fuel efficiency. 

Some literature can be found where the optimization problem is formulated to minimize the fuel consumption rate or fuel cost directly. In \cite{hernandez-aramburo_fuel_2005}, four power-sharing schemes are presented to establish a unit commitment strategy to minimize fuel costs. However, instead of considering the generator efficiency function, the authors utilized a typical fuel consumption characteristic curve, which largely depends on the capacity and model of a generator. 
The authors in \cite{lundh_estimation_2016} implement a recursive method to estimate a second-order polynomial model of specific fuel consumption. The model is later used to determine the optimal load distribution between the different generators. A similar approach is found in \cite{kanellos_optimal_2014}, where dynamic programming is used to solve the formulated problem. However, both optimization models are suitable for simple power networks since they don't include AC power flow or energy storage models. 
A minimum hourly fuel consumption curve interpolated by a quadratic equation is used in \cite{barsali_control_2004} to minimize the fuel consumption of hybrid electric vehicles. This method can not be extended for transmission or distribution systems since most power system constraints are not included in the model. 
In \cite{accetta_energy_2019}, the authors utilize a 3D map of brake-specific fuel consumption (BSFC) in terms of the rotating speed of the drive and generated mechanical torque to determine the minimum point of diesel engine (DE) fuel consumption. Authors in \cite{satpathi_modeling_2017} also apply a similar method where a speed vs. power curve of the diesel engine (DE) is exploited to achieve minimum fuel consumption conditions. Nonetheless, these methods have two major drawbacks: they are implemented particularly for the DC ship systems, and the strong nonlinearity of the BSFC curve increases the complexity of determining the optimal operation. 

Reduction in fuel cost can also be achieved indirectly through economic dispatch (ED) optimization problems \cite{chowdhury_review_1990}, \cite{kunya_review_2023}. In the ED problem, a polynomial objective function (generally in quadratic form) representing the cost of the generator dispatch is optimized to supply the demand most economically. However, since the objective function of the ED problem does not represent the fuel consumption rate, these formulations are unable to accomplish the highest system efficiency. They are also inconvenient for long-term planning of fuel usage and impractical to apply to systems where achieving optimal fuel consumption rate is the main goal. 

Although the strategies discussed above for improving overall system performance and reducing fuel cost have their own merits, the key to minimizing fuel consumption lies in introducing a function that accurately represents the fuel consumption rate. By directly representing fuel consumption in the optimization process, researchers can effectively tackle the core challenge of reducing fuel usage and achieving greater energy efficiency in the system. Therefore, future studies may benefit from exploring methodologies considering this crucial aspect when addressing system efficiency optimization.

This paper presents a novel approach by introducing a unique sum-of-ratios objective function, which directly represents the fuel consumption rate. Unlike conventional methods that optimize the polynomial generator efficiency function or cost of operation, this sum-of-ratios formulation offers a practical and more efficient solution. However, solving multiple ratio optimization problems, known as Fractional Programming (FP), has been proven to be NP-hard \cite{freund_solving_nodate}. The convergence of these problems with established nonlinear optimization methods can take an extensive amount of time. In addition, when the sum-of-ratios problem involves more than 20 ratios, the current approaches struggle to find a solution within a reasonable timeframe \cite{kuno_branch-and-bound_nodate},\cite{carlsson_linear_2013}.
Due to these challenges, directly solving the sum-of-ratios optimization problem for energy management systems (EMS) is not feasible. Especially this method will be inapplicable for large-scale systems with a high number of variables. It necessitates the development of an efficient reformulation and solution technique to address multiple ratio problems effectively. Finding innovative approaches to tackle these difficulties will be crucial in making the proposed sum-of-ratios method a practical and scalable solution for optimizing fuel consumption in real-world energy systems.

A substantial body of literature exists on Fractional Programming (FP), but the emphasis has primarily been on single-ratio problems. Among the renowned reformulation techniques, the \emph{Charnes-Cooper Transform}, \cite{charnes_programming_1962} \cite{schaible_parameter-free_1974} is notable for proposing an algorithm to solve single ratio linear FP problems by introducing two new variables and converting the fractional problem into a linear problem. Another classical technique, \emph{Dinkelbach's Transform} \cite{dinkelbach_nonlinear_1967}, reformulates the single ratio problem using a new auxiliary variable updated iteratively until convergence is achieved.
\cite{toloo_equivalent_2021} presents a formulation of a linear problem equivalent to a single ratio linear FP problem where some duality properties are used to prove the equivalence. For quadratic FP problems, where both the numerator and denominator are quadratic functions, a new method called the \emph{decomposition fractional separable method} is proposed in \cite{jayalakshmi_new_nodate} using linear programming techniques. 
An alternative approach to solving single-ratio quadratic FP is outlined in \cite{department_of_mathematical_engineering_yildiz_technical_university_istanbul_turkey_novel_2018}, employing Taylor series expansion for effective reformulation.

The literature discussed in references \cite{charnes_programming_1962}-\cite{department_of_mathematical_engineering_yildiz_technical_university_istanbul_turkey_novel_2018} primarily focused on solving single-ratio FP problems and cannot be directly extended to handle multi-ratio problems. Addressing multiple ratio problems, as encountered in the sum-of-ratios function, remains a challenge that requires innovative and efficient reformulation and solution techniques. However, in \cite{almogy_class_1971}, the authors proposed an extension of Dinkelbach's Transform specifically tailored to address multi-ratio FP problems. Nonetheless, this method was later refuted by Falk and Palocsay \cite{floudas_optimizing_1991}, who demonstrated its limitations through a numerical example.

To find the globally optimal solution for the sum-of-ratios problem, \cite{jong_practical_nodate} introduced a practical method that involves solving a sequence of convex programming problems. In \cite{tsai_global_2005}, a convexification strategy was employed to decompose fractional terms into convex and concave components. Then, a piecewise linearization technique was applied to approximate the concave terms effectively.
Additionally, \cite{shen_fractional_2018} proposed a quadratic transform to tackle concave-convex multiple ratio minimization problems. In the case of generalized convex multiplicative functions, a reformulation technique was presented in \cite{konno_global_1994}, where the main problem was reformulated as a concave minimization problem with 2p variables. This reformulation technique could also be applied to sum-of-ratios FP problems if the multiplicative terms were replaced with a convex over a concave function \cite{schaible_fractional_2003}. 

In our study, the objective function is defined as a sum-of-ratios minimization problem with non-negative-convex numerator and positive-concave denominator terms. Due to this specific form of the problem, an appropriate algorithm should be selected to solve the formulated multiple ratios FP problem effectively. As a result, the reformulation technique presented in \cite{konno_global_1994} is adopted in this paper. By leveraging this reformulation technique, the complexities of the sum-of-ratios problem can be effectively addressed, and an optimized solution can be found with a feasible convergence time. 

The contributions of this paper are the following: 
\begin{itemize}
\item A novel fractional objective function is introduced in this literature, which directly represents the fuel consumption rate of the generators. Unlike typical system efficiency optimization problems that use the efficiency function or the operating cost function as the objective, this unique formulation directly accounts for fuel consumption. This approach proves to be more efficient and practical compared to previous studies since it directly targets the core issue of minimizing fuel usage and improving overall system efficiency. 

\item To address the optimization problem in this study, the sum-of-ratios fractional programming (FP) algorithm is employed. to the best of our knowledge, this literature represents the first application of the FP method to solve the optimization problem for EMS efficiently. The reformulation technique with FP can also be applied to different power or communication system research where sum-of-ratios functions are used.  

\item The successful application of the FP algorithm, combined with the convex relaxation of nonlinear constraints, demonstrates that the proposed model is suitable for handling large-scale systems. This capability is exemplified through the model's effective implementation on the IEEE 118-bus system. By demonstrating its applicability to such a complex and extensive system, the paper establishes the scalability of the proposed approach for real-world energy management scenarios.
\end{itemize}

The remainder of the paper is organized as follows: the fuel efficiency problem formulation with sum-of-ratios objective function and its convex reformulation are presented in Section \ref{System_Efficiency:section2}. In Section \ref{System_Efficiency:section3}, the solution algorithm for the reformulated problem is described. The results for the notional MVAC ship system, IEEE 30-bus, and IEEE 118-bus system, accompanied with the performance comparisons, are demonstrated in Section \ref{System_Efficiency:section4}. Finally, Section \ref{System_Efficiency:section5} represents the conclusion and future work.

\section{Problem Formulation} \label{System_Efficiency:section2}

\subsection{Optimization Model for System Efficiency} \label{System_Efficiency:section2.1}
This section presents a unique sum-of-ratios objective function for the optimization problem that directly represents the fuel consumption of the generators. The objective of the minimization problem is to minimize the fuel consumption rate of the generators over the planning horizon, which will maximize the system efficiency. The fuel consumed by a generator can be expressed by taking into account the generator's efficiency and the output power it produces. The objective function is the following:
\begin{equation}
	\text{minimize } f = \frac{1}{\alpha}\sum_{t=0}^{T} \sum_{i\in N_g} \frac{P_{i,g}^t}{\eta_{i,g}} \Delta t. 
	\label{eq:System_Efficiency:Eq1}
\end{equation}
where, generator efficiency, $\eta = a_i p_{i,g}^2 + b_i p_{i,g} + c_i$ with $a$, $b$, and $c$ are generator specific constants, $\alpha$ is the fuel energy density (MWh/L), 
$p_{i,g}$ is the per-unit output of the i-th generator: $p_{i,g} = {P_{i,g}}/{P_{i,b}}$; 
$P_{i,g}^t$ is the generator output power at time t, and $P_{i,b}$ is the base power of i-th generator.
$N_{g}$ = total number of generators, $T$ = planning horizon, $\triangle t$ = each time period.\

The following active and reactive power balance constraints are associated with the system: 
\begin{equation}
    P_{i,inj} = P_{i,g} - P_{i,l}
    \label{eq:System_Efficiency:Eq2}
\end{equation}
\begin{equation}
    Q_{i,inj} = Q_{i,g} - Q_{i,l}
    \label{eq:System_Efficiency:Eq3}
\end{equation}

where 
\begin{equation}
    P_{i,inj} = \sum_{k\in B} V_i V_k (G_{ik} cos\theta_{ik} + B_{ik} sin\theta_{ik})
    \label{eq:System_Efficiency:Eq4}
\end{equation}

\begin{equation}
    Q_{i,inj} = \sum_{k\in B} V_i V_k (G_{ik} sin\theta_{ik} + B_{ik} cos\theta_{ik})
    \label{eq:System_Efficiency:Eq5}
\end{equation}

where (\ref{eq:System_Efficiency:Eq2}) and (\ref{eq:System_Efficiency:Eq5}) are the AC power flow constraints for the system. 
The following constraints should be included in the problem formulation to maintain the operational limits of the system: 

\begin{equation}
    P_{i,g}^{min} \leq P_{i,g}^t \leq P_{i,g}^{max},
    \label{eq:System_Efficiency:Eq6}
\end{equation}

\begin{equation}
    Q_{i,g}^{min} \leq Q_{i,g}^t \leq Q_{i,g}^{max},
    \label{eq:System_Efficiency:Eq7}
\end{equation}

\begin{equation}
    V_{i}^{min} \leq V_{i}^t \leq V_{i}^{max},
    \label{eq:System_Efficiency:Eq8}
\end{equation}

\begin{equation}
    \theta_{i}^{min} \leq \theta_{i}^t \leq \theta_{i}^{max},
    \label{eq:System_Efficiency:Eq9}
\end{equation}

\begin{equation}
    -P_{ik} \leq P_{ik}^t \leq P_{ik},
    \label{eq:System_Efficiency:Eq10}
\end{equation}

\begin{equation}
    -Q_{ik} \leq Q_{ik}^t \leq Q_{ik},
    \label{eq:System_Efficiency:Eq10.1}
\end{equation}

\begin{equation}
    -DR_{i} \leq P_{i,g}^{t+1} - P_{i,g}^{t} \leq UR_{i}.
    \label{eq:System_Efficiency:Eq11}
\end{equation}

where (\ref{eq:System_Efficiency:Eq6}) and (\ref{eq:System_Efficiency:Eq7}) represent the generators' real and reactive power generation limits, (\ref{eq:System_Efficiency:Eq8}) and (\ref{eq:System_Efficiency:Eq9}) are the voltage and voltage angle limits, (\ref{eq:System_Efficiency:Eq10}) and (\ref{eq:System_Efficiency:Eq10.1}) are the line limits for real and reactive power, and (\ref{eq:System_Efficiency:Eq11}) is the ramp rate limit. 

The energy storage system (ESS) plays a vital role in minimizing the fuel consumption by the generators. The following Energy Storage System (ESS) constraints are included in the optimization problem:
\begin{equation}
    E_{i,b}^{t} = E_{i,b}^{t-1} - \eta_b P_{i,b}^{r,t} \triangle t
    \label{eq:System_Efficiency:Eq12}
\end{equation}

\begin{equation}
    -P_{i}^{C,max} \leq P_{i,b}^{r,t} \leq P_{i}^{D,max},
    \label{eq:System_Efficiency:Eq13}
\end{equation}

\begin{equation}
    \sum_{t=0}^{t} P_{i,b}^{r,t} = 0,
    \label{eq:System_Efficiency:Eq14}
\end{equation}

\begin{equation}
    SOC_{i}^{min} \leq SOC_{i}^{t} \leq SOC_{i}^{max},
    \label{eq:System_Efficiency:Eq15}
\end{equation}

where (\ref{eq:System_Efficiency:Eq12}) indicates the energy conservation constraint, (\ref{eq:System_Efficiency:Eq13}) is the limit for charging or discharging rate, and (\ref{eq:System_Efficiency:Eq15}) is the state of charge (SOC) limit of the ESS. For the ESS, although the SOC can vary from 0 to 1 (0\% to 100\%), fully discharging can damage the battery permanently and shorten the life cycle of the battery \cite{noauthor_battery_nodate}. In this paper, the minimum SOC is selected as 0.2 (20\%). Eq.(\ref{eq:System_Efficiency:Eq14}) ensures that the sum of the total charging and discharging power over a planning period will be zero, which helps the system to recharge the battery before the next planning cycle.

\subsection{Reformulation of the Objective Function} \label{System_Efficiency:section2.2}
The objective function (\ref{eq:System_Efficiency:Eq1}) for the formulated optimization problem is a highly nonlinear sum-of-ratios case. In this section, the objective function is reformulated to a convex function using a technique used for solving generalized convex multiplicative problems. Later, the convex minimization method is utilized to solve the reformulated problem iteratively.

The general convex multiplicative minimization problem has the following structure: 

\begin{equation}
    \text{minimize } h(x) + \sum_{i=1}^{p} f_i(x)g_i(x)
    \label{eq:System_Efficiency:Eq16}
\end{equation}

\centerline{\text{subject to } $x \in X$}

where $h, f_i(x)$ and $g_i(x)$ for all $i$ are convex functions and $X\subset R^n$ is a convex set.
If $h(x)=0$, $f_i(x)=A_i(x)$, and $g_i(x)=1/B_i(x)$, (\ref{eq:System_Efficiency:Eq16}) will be in the following form: 

\begin{equation}
    \text{minimize } H(x) = \sum_{i=1}^{m} \frac{A_i(x)}{B_i(x)}
    \label{eq:System_Efficiency:Eq17}
\end{equation}

\centerline{\text{subject to } $x \in X$}
which is a sum-of-ratios problem, where $A_i(x)$ are non-negative, convex and $B_i(x)$ are positive, concave functions for all $i$.

The authors in \cite{konno_global_1994} defined the following problem by introducing 2m auxiliary variables $\zeta_i$ and $\beta_i$, where $i=1,2,3,....m$: 
\begin{equation}
    \text{minimize } F(x, \zeta, \beta) = \frac{1}{2} \sum_{i=1}^{m} [\zeta_i (A_i(x))^2 + \beta_i (B_i(x))^2]
    \label{eq:System_Efficiency:Eq18}
\end{equation}

\centerline{\text{subject to } $x \in X$}
\centerline{$\zeta_i \beta_i \ge 1$}
\centerline{$(\zeta, \beta) > 0$}

where $\zeta = (\zeta_1, \zeta_2, ....\zeta_m)$ and $\beta = (\beta_1, \beta_2, ..... \beta_m)$. It can be proved that, if $(x^*, \zeta ^*, \beta ^*)$ is an optimal solution of (\ref{eq:System_Efficiency:Eq18}), then $x^*$ will be an optimal solution of (\ref{eq:System_Efficiency:Eq17}) and $H(x^*) = F(x^*, \zeta ^*, \beta ^*)$ \cite{konno_global_1994}.

As a result, the optimization problem in section \ref{System_Efficiency:section2}(A) can be written as the following problem: 

\begin{equation}
    \text{minimize } f(P, \zeta, \beta) = \frac{1}{2} [\sum_{i=1}^{m} (\zeta_i (P_{i,g}^t)^2 + \beta_i \eta_{i,g}^2)] \Delta t
    \label{eq:System_Efficiency:Eq19}
\end{equation}
subject to\\
\centerline{(\ref{eq:System_Efficiency:Eq2})-(\ref{eq:System_Efficiency:Eq15})},

\begin{equation}
    \zeta_i \beta_i \ge 1,
    \label{eq:System_Efficiency:Eq20}
\end{equation}

\begin{equation}
    (\zeta, \beta) > 0,
    \label{eq:System_Efficiency:Eq21}
\end{equation}

\subsection{Convex Relaxation Technique} \label{System_Efficiency:section2.3}

The convex optimization methods can only be applied to problems where the objective function and all constraints are finite and convex. Although the reformulated problem in section \ref{System_Efficiency:section2}(A) has a convex objective function for a fixed set of $\eta_i$ and $\beta_i$, several constraints still have nonlinearity. In this section, the nonlinear power flow constraints (\ref{eq:System_Efficiency:Eq2}), (\ref{eq:System_Efficiency:Eq3}), and line flow constraints (\ref{eq:System_Efficiency:Eq10}),(\ref{eq:System_Efficiency:Eq10.1}) are replaced with the following linear and quadratic constraints: 

\begin{equation}
    P_{i,\mathrm{inj}} = \sqrt{2}u_i G_{ii} + \sum_{k \in B} \left(G_{ik}W_{R_{ik}} + B_{ik}W_{I_{ik}}\right),
    \label{eq:System_Efficiency:Eq22}
\end{equation}

\begin{equation}
    Q_{i,\mathrm{inj}} = -\sqrt{2}u_i B_{ii} + \sum_{k \in B} \left(G_{ik}W_{I_{ik}} - B_{ik}W_{R_{ik}}\right),
    \label{eq:System_Efficiency:Eq23}
\end{equation}

\begin{equation}
    P_{ik} = \sqrt{2}u_i G_{ik} - \left(G_{ik}W_{R_{ik}} + B_{ik}W_{I_{ik}}\right),
    \label{eq:System_Efficiency:Eq24}
\end{equation}

\begin{equation}
    Q_{ik} = -\sqrt{2}u_i B_{ik} + \left(B_{ik}W_{R_{ik}} - G_{ik}W_{I_{ik}}\right).
    \label{eq:System_Efficiency:Eq25}
\end{equation}

\begin{equation}
    W_{R_{ik}}^2 + W_{I_{ik}}^2 \leq 2u_i u_k.
    \label{eq:System_Efficiency:Eq26}
\end{equation}

\begin{equation}
    \theta_{ik} = \tan^{-1} \left( \frac{W_{I_{ik}}}{W_{R_{ik}}} \right).
    \label{eq:System_Efficiency:Eq27}
\end{equation}

Here, (\ref{eq:System_Efficiency:Eq22}), (\ref{eq:System_Efficiency:Eq23}) are the linear real and reactive power flow equations and (\ref{eq:System_Efficiency:Eq24}), (\ref{eq:System_Efficiency:Eq25}) are linear line flow equations. The relationship between the convex variables $u_i, W_{I_{ik}}$, and $W_{R_{ik}}$ are defined by the equations (\ref{eq:System_Efficiency:Eq26}) and (\ref{eq:System_Efficiency:Eq27}). Since (\ref{eq:System_Efficiency:Eq27}) is still nonlinear, Taylor series expansion can be used to linearize the equation: 

\begin{equation}
\begin{split}
    \tan^{-1}\frac{W_{I_{ik}}^{(q)}}{W_{R_{ik}}^{(q)}}=\theta_{ik}+\frac{W_{I_{ik}}^{(q)}}{W_{R_{ik}}^{(q)2}+W_{I_{ik}}^{(q)2}} W_{R_{ik}}- \\ \frac{W_{R_{ik}}^{(q)}}{W_{R_{ik}}^{(q)2}+W_{I_{ik}}^{(q)2}}W_{i_{ik}}.
    \label{eq:System_Efficiency:Eq28}
\end{split}
\end{equation}

where, the higher order terms are neglected, and $(W_{I_{ik}}^{(q)}, W_{R_{ik}}^{(q)})$ are the initial estimation. A detailed description of this technique can be found in \cite{anam_system_2023}.\\
The final problem formulation will be as follows: \\

\centerline{\text{minimize $f(x, \zeta, \beta)$}}

subject to\\
\centerline{(\ref{eq:System_Efficiency:Eq2}), (\ref{eq:System_Efficiency:Eq3}), (\ref{eq:System_Efficiency:Eq6})-(\ref{eq:System_Efficiency:Eq15}), (\ref{eq:System_Efficiency:Eq20}), (\ref{eq:System_Efficiency:Eq21}), (\ref{eq:System_Efficiency:Eq26}) and (\ref{eq:System_Efficiency:Eq28})}.

where, $P_{i,inj}, Q_{i,inj}, P_{ik}$, and $Q_{ik}$ are defined by (\ref{eq:System_Efficiency:Eq22}), (\ref{eq:System_Efficiency:Eq23}), (\ref{eq:System_Efficiency:Eq24}), and (\ref{eq:System_Efficiency:Eq25}), respectively.

\section{Solution Technique} \label{System_Efficiency:section3}

In this section, an iterative method is described to solve the reformulated problem. 
For a fixed set of $(\zeta, \beta)$, let us consider the following sub-problem of (\ref{eq:System_Efficiency:Eq18}):

\begin{equation}
    \text{minimize } F(x; \zeta, \beta) = \frac{1}{2} \sum_{i=1}^{m} [\zeta_i (A_i(x))^2 + \beta_i (B_i(x))^2]
    \label{eq:System_Efficiency:Eq29}
\end{equation}

Equation (\ref{eq:System_Efficiency:Eq29}) can be solved using any standard convex optimization technique. If the optimal solution for (\ref{eq:System_Efficiency:Eq29}) is $x^*(\zeta, \beta)$ , then for a fixed set of $x^*$ (\ref{eq:System_Efficiency:Eq17}) is reduced to the following problem of 2m variables $(\zeta, \beta)$:

\begin{equation}
    \text{minimize } F_{aux}(\zeta, \beta) = \frac{1}{2} \sum_{i=1}^{m} [\zeta_i (A_i(x^*))^2 + \beta_i (B_i(x^*))^2]
    \label{eq:System_Efficiency:Eq30}
\end{equation}

subject to,

\centerline{$\zeta_i \beta_i \ge 1$}
\centerline{$(\zeta, \beta) > 0$}

where, $F_{aux}$ denotes the auxiliary problem of $F$. Equation (\ref{eq:System_Efficiency:Eq29}) and (\ref{eq:System_Efficiency:Eq30}) are solved iteratively until convergence is achieved. 

The following algorithm is used to solve the fractional optimization problem: 

\begin{algorithm}[H]
\caption{Fractional Programming Solution Steps}\label{alg:alg1}
\begin{algorithmic}
\STATE 
\STATE \hspace{0.1cm}$ \textbf{Step 0: }$ \textbf{Set i = 0}.
\STATE \hspace{0.1cm}$ \textbf{Step 1: }$ \textbf{Find an optimal solution, $x^*$ for (\ref{eq:System_Efficiency:Eq29}) for a fixed $(\zeta, \beta)$}.
\STATE \hspace{0.1cm}$ \textbf{Step 2: }$ \textbf{Find an optimal solution $(\zeta_i, \beta_i)$ for (\ref{eq:System_Efficiency:Eq30}) using convex minimization method}.
\STATE \hspace{0.1cm}$ \textbf{Step 4: }$ \textbf{Update the feasible region of (\ref{eq:System_Efficiency:Eq30}) by including cutting function $l_i$, where $l_i = 2-\beta_i \sqrt\frac{\zeta_i}{\beta_i} - \zeta_i\sqrt\frac{\beta_i}{\zeta_i} \le0$}
\STATE \hspace{0.1cm}$ \textbf{Step 5: }$ \textbf{Let $i=i+1$ and return to step 0. Continue until converges}.
\end{algorithmic}
\label{alg1}
\end{algorithm}

In this paper, the MOSEK optimization toolbox is used to solve the formulated problem \cite{noauthor_mosek_nodate}. The optimization problem is transformed into the conic quadratic format to fit into MOSEK. MOSEK supports two types of quadratic cones: 

\begin{itemize}
  \item General quadratic cones:\\ 
    \centerline{$Q^n = [x \in R^n: x_0\ge\sqrt{\sum_{j=1}^{n-1} x_j^2}]$}\\
    
  \item Rotated Quadratic cones:\\
  \centerline{$Q_r^n = [x \in R^n: 2x_0 x_1\ge\sum_{j=2}^{n-1} x_j^2], x_0 \ge 0, x_1 \ge 0$}
\end{itemize}

All the quadratic parts of the reformulated minimization problem in section \ref{System_Efficiency:section2}(C) are replaced with rotated quadratic cones and corresponding linear equations. As a result, the transformed problem can be solved using the MOSEK solver. 

\section{Case Studies} \label{System_Efficiency:section4}

\begin{figure}[t]
\centering
\includegraphics[scale=0.11]{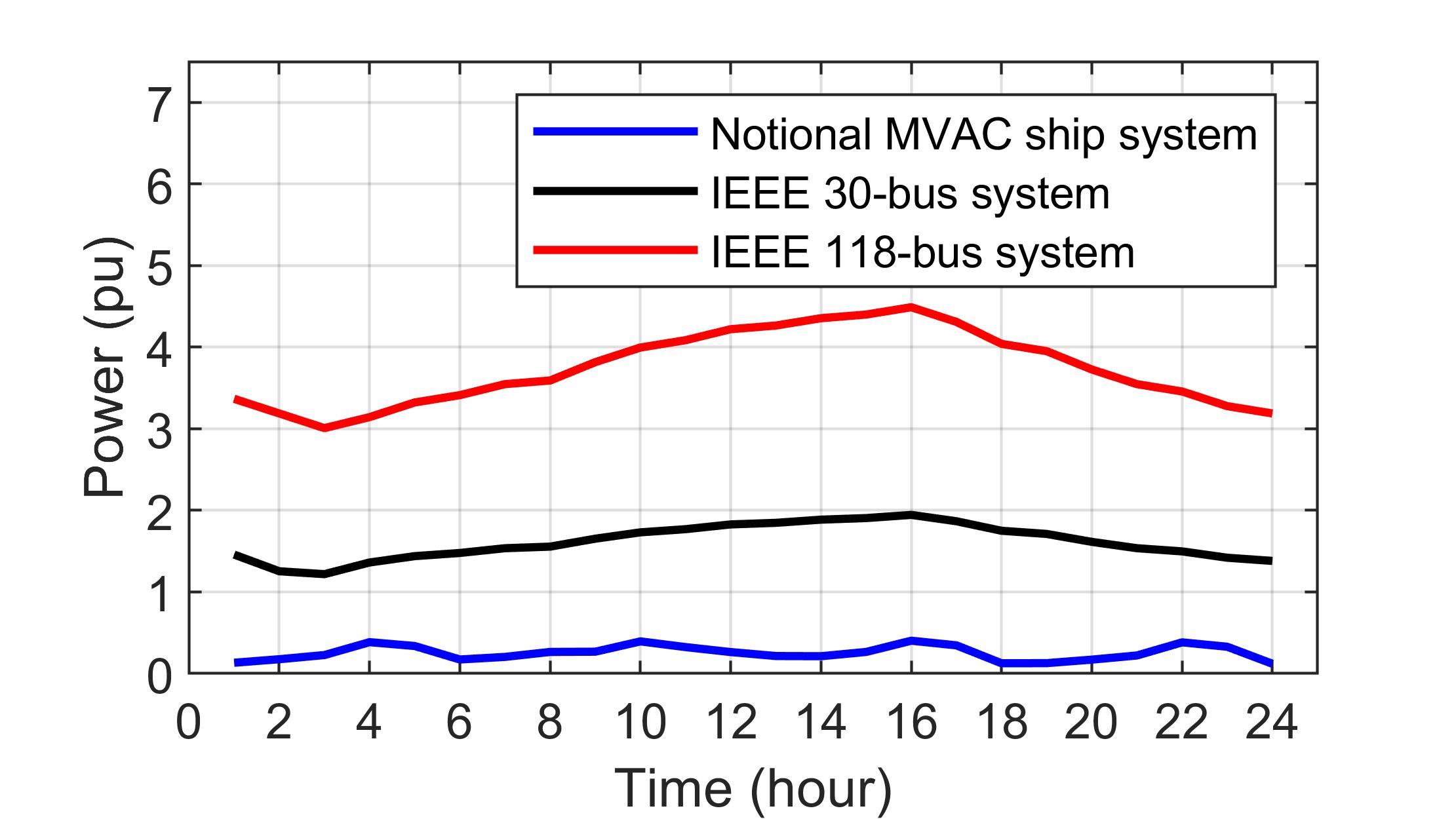}
\caption{24-hours load profile of different systems}
\label{fig:load_profile}
\end{figure}

In this literature, the proposed system efficiency model is tested with a notional 12-bus MVAC ship system, a modified IEEE 30-bus, and an IEEE 118-bus system. Each system consists of multiple generators, distribution lines, ESS, and loads at different buses. The load data is generated using the demand pattern of Real-time Dashboard, NYISO \cite{noauthor_real-time_nodate}. The load profile for different systems can be observed in Fig.\ref{fig:load_profile}. This paper considers a 24-hour time horizon while solving the optimization problem, where each time step is 1 hour. However, any time horizon and length of time step can be selected based on the system requirement.

\subsection{Notional 12-Bus MVAC Ship System} \label{System_Efficiency:section4.1}
\begin{figure}[t]
\centering
\includegraphics[scale=0.27]{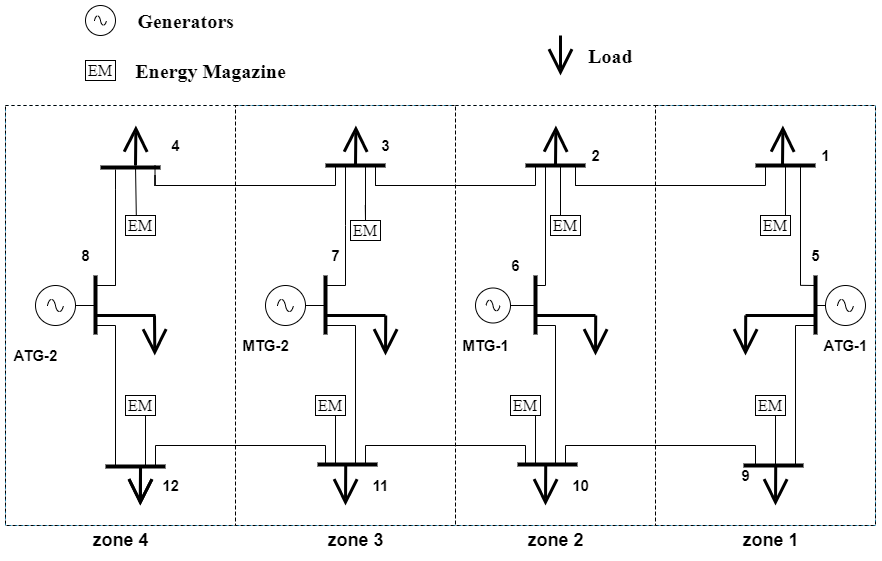}
\caption{A notional MVAC ship system model}
\label{fig:MVAC_model}
\end{figure}

The notional four-zone 12 bus MVAC ship system \cite{noauthor_esrdc_nodate} (shown in Fig. \ref{fig:MVAC_model}) consists of two main gas turbine generators (MTG) and two auxiliary gas turbine generators (ATG). The generator parameters can be observed from TABLE \ref{tab:MVAC_gen_data}. The system also has 8 ESS (each zone has 2 ESS) and multiple loads, including 2 propulsion motor modules (PMM) at buses 6 and 7 and AC load centers (ACLC) at buses 1,2,3,4,9,10,11, and 12. The ESS data are shown in TABLE \ref{tab:ESS_data}.

\begin{table}[b]
	\renewcommand{\arraystretch}{1.3}
	\caption{Generator Data for Notional MVAC Ship System}
	\label{tab:MVAC_gen_data}
	\centering
	\begin{tabular}
		{
			>{\centering\arraybackslash}m{0.08\textwidth}
			>{\centering\arraybackslash}m{0.08\textwidth}
			>{\centering\arraybackslash}m{0.14\textwidth}
			>{\centering\arraybackslash}m{0.08\textwidth}
		}
		\hline
		Types 
		& Capacity (MW)
		& Efficiency curve coefficient
		& Number of Units 
		\\
		\hline
		ATG & 4.7 & $a_i = -.133, b_i = .311, c_i = .174$ & 2 \\
        MTG & 35 & $a_i = -.133, b_i = .311, c_i = .204$ & 2 \\
		\hline\hline
        \end{tabular}
\end{table}

\begin{figure}
\centering
\includegraphics[scale=0.48]{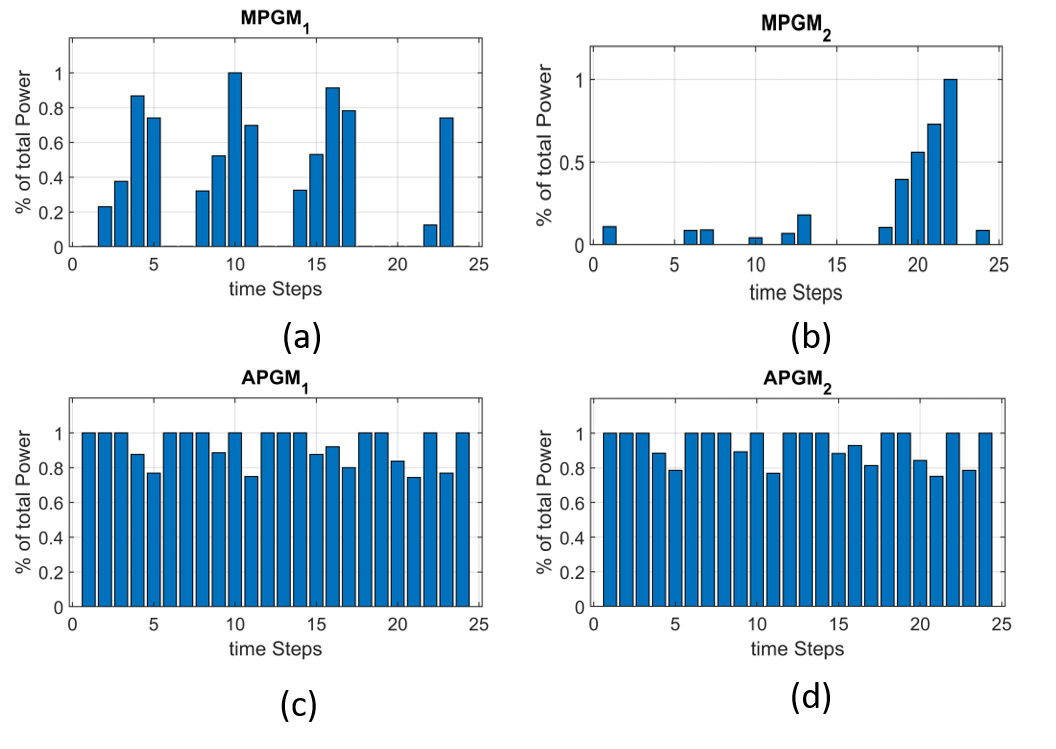}
\caption{Generator output for notional MVAC ship system: (a) MPGM 1, (b) MPGM 2, (c) APGM 1, (d) APGM 2.}
\label{fig:MVAC_gen}
\end{figure}

\begin{figure}
\centering
\includegraphics[scale=0.48]{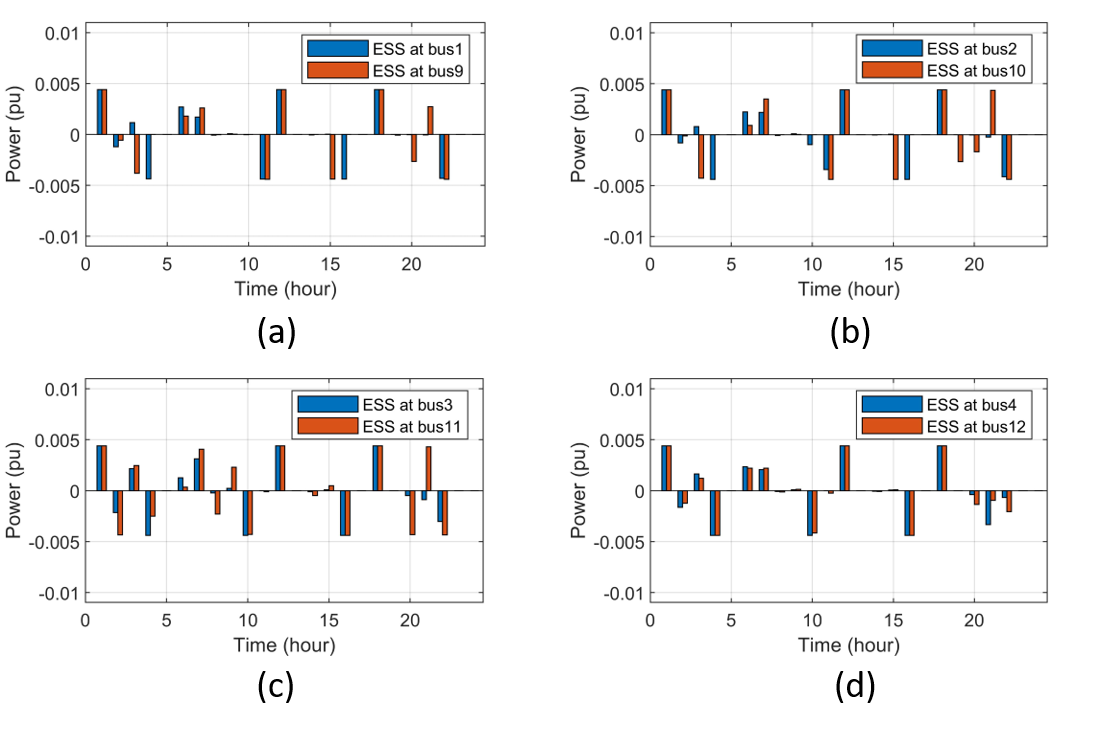}
\caption{ESS charging and discharging schedule for notional MVAC ship system; (a) Zone 1, (b) Zone 2, (c) Zone 3, (d) Zone 4.}
\label{fig:MVAC_ess}
\end{figure}

\begin{figure}
\centering
\includegraphics[scale=0.08]{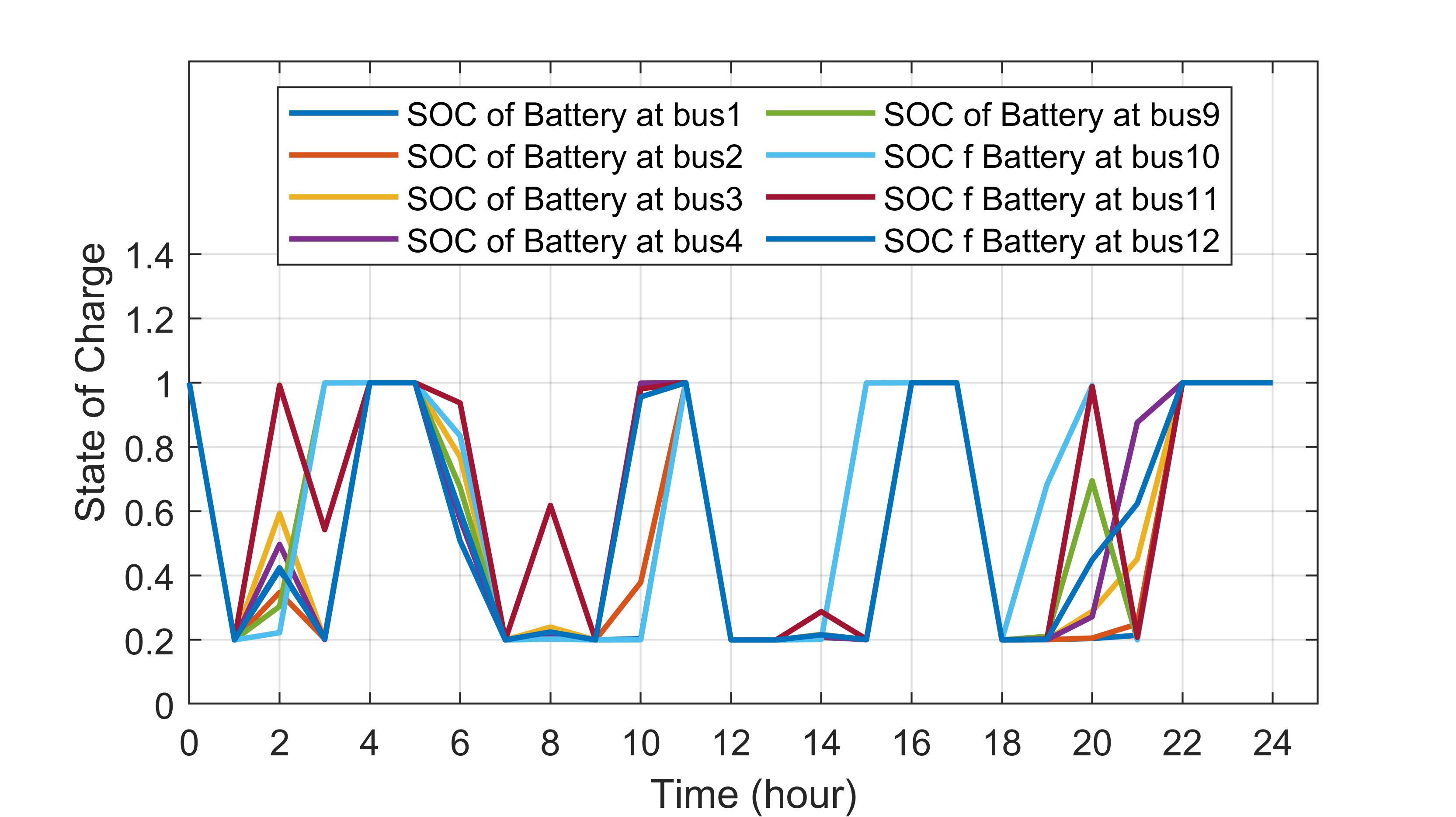}
\caption{SOC of ESS for notional MVAC ship system}
\label{fig:MVAC_soc}
\end{figure}

The proposed optimization model was run for the notional MVAC ship system; the simulation time was 34.414s, taking 178 iterations to converge. The output power generation can be observed in Fig. \ref{fig:MVAC_gen}. The system is likely to dispatch the ATGs more than the MTGs to improve the overall efficiency since the capacities of the ATGs are lower than the MTGs.
The ESS charging and discharging schedule and the SOC of the ESS are shown in Figs. \ref {fig:MVAC_ess} and \ref{fig:MVAC_soc}, respectively. 

\begin{table}[b]
	\renewcommand{\arraystretch}{1.3}
	\caption{ESS data for notional MVAC ship system, IEEE 30-bus, and IEEE 118-bus system}
	\label{tab:ESS_data}
	\centering
	\begin{tabular}
		{
			>{\centering\arraybackslash}m{0.1\textwidth}
			>{\centering\arraybackslash}m{0.05\textwidth}
			>{\centering\arraybackslash}m{0.05\textwidth}
			>{\centering\arraybackslash}m{0.1\textwidth}
			>{\centering\arraybackslash}m{0.08\textwidth}
		}
		\hline
		System 
		& Number of ESS
		& Capacity (MWh) 
		& Maximum Charging/discharging rate (MW/h) 
		& Miniumum SOC (\%) \\
		\hline
		MVAC ship system & 8 & 2.2 & 10 & 20\\
        IEEE 30-bus system & 6 & 15 & 5 & 20\\
        IEEE 118-bus system & 16 & 15 & 5 & 20\\
		\hline\hline
        \end{tabular}
\end{table}

\subsection{IEEE 30-bus System} \label{System_Efficiency:section4.2}

The IEEE 30-bus system \cite{noauthor_ieee_nodate} has 6 generator buses, 16 load buses, and 42 transmission lines. In addition, the system is modified by including six energy storage systems (ESS) at buses 5, 11, 15, 19, 23 and 27. 

\begin{figure}[t]
\centering
\includegraphics[scale=0.55]{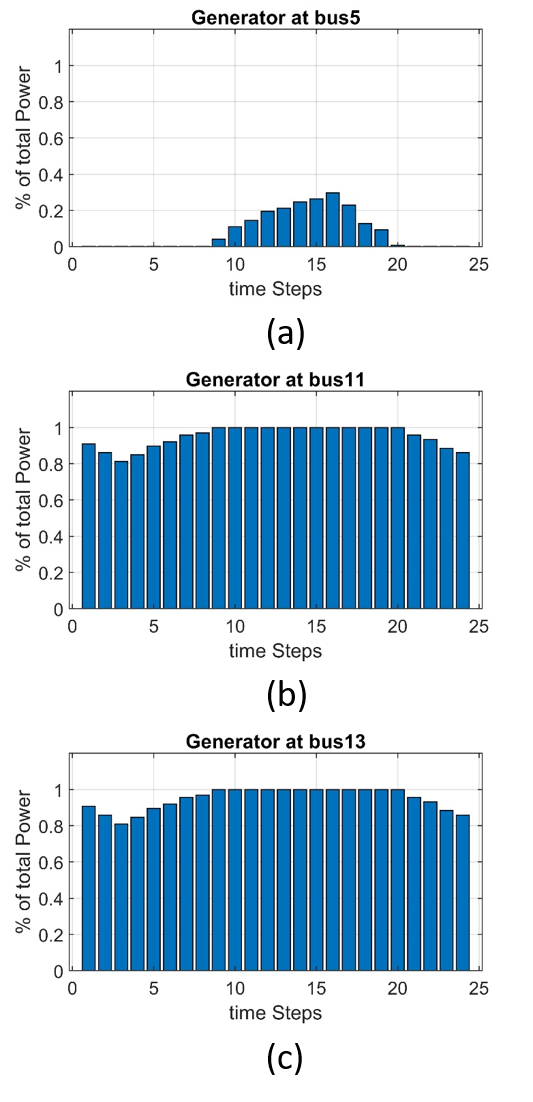}
\caption{Generator output for IEEE 30-bus system at (a) bus 5, (b) bus 11, (c) bus 13.}
\label{fig:gen_30}
\end{figure}

\begin{figure}
\centering
\includegraphics[scale=0.08]{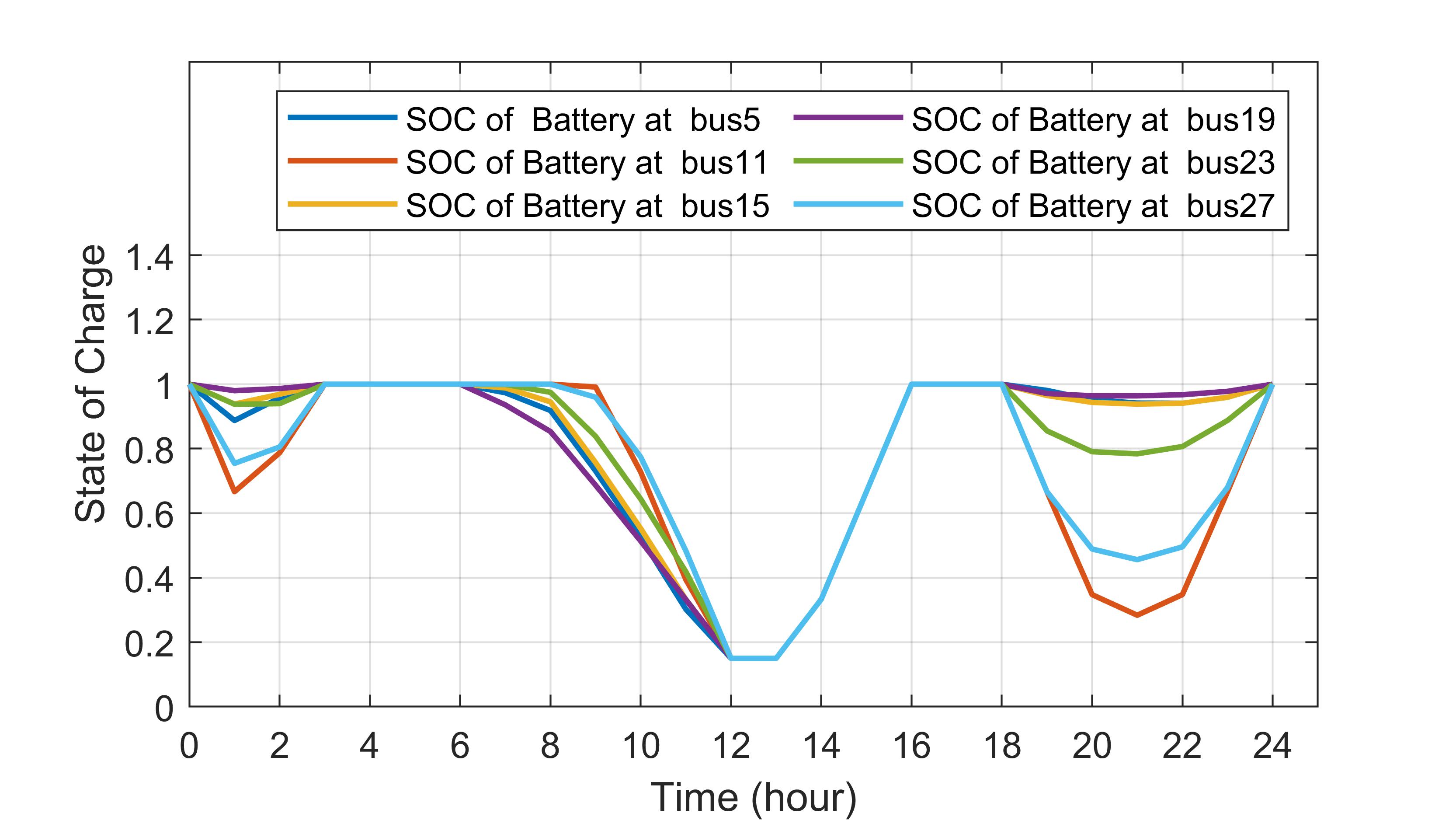}
\caption{State of charge (soc) for IEEE 30-bus system}
\label{fig:soc_30}
\end{figure}

The proposed model was then tested for the IEEE 30-bus system, and the observed convergence time was 91.728s for 476 iterations. The generator output schedule is shown in Fig. \ref{fig:gen_30}. Although the system has six generators, only three produced output during the simulation. Generator graphs with zero output are not included in the figure. 
The state of charge (SOC) of the ESS is shown in Fig. \ref {fig:soc_30}.

\subsection{IEEE 118-bus System} \label{System_Efficiency:section4.2}

The IEEE 118-bus system \cite{noauthor_ieee_nodate-1} has 21 generator buses, 113 load buses, and 179 transmission lines. In addition, the system is modified by including 16 energy storage systems (ESS). 

\begin{figure}[t]
\centering
\includegraphics[scale=0.3]{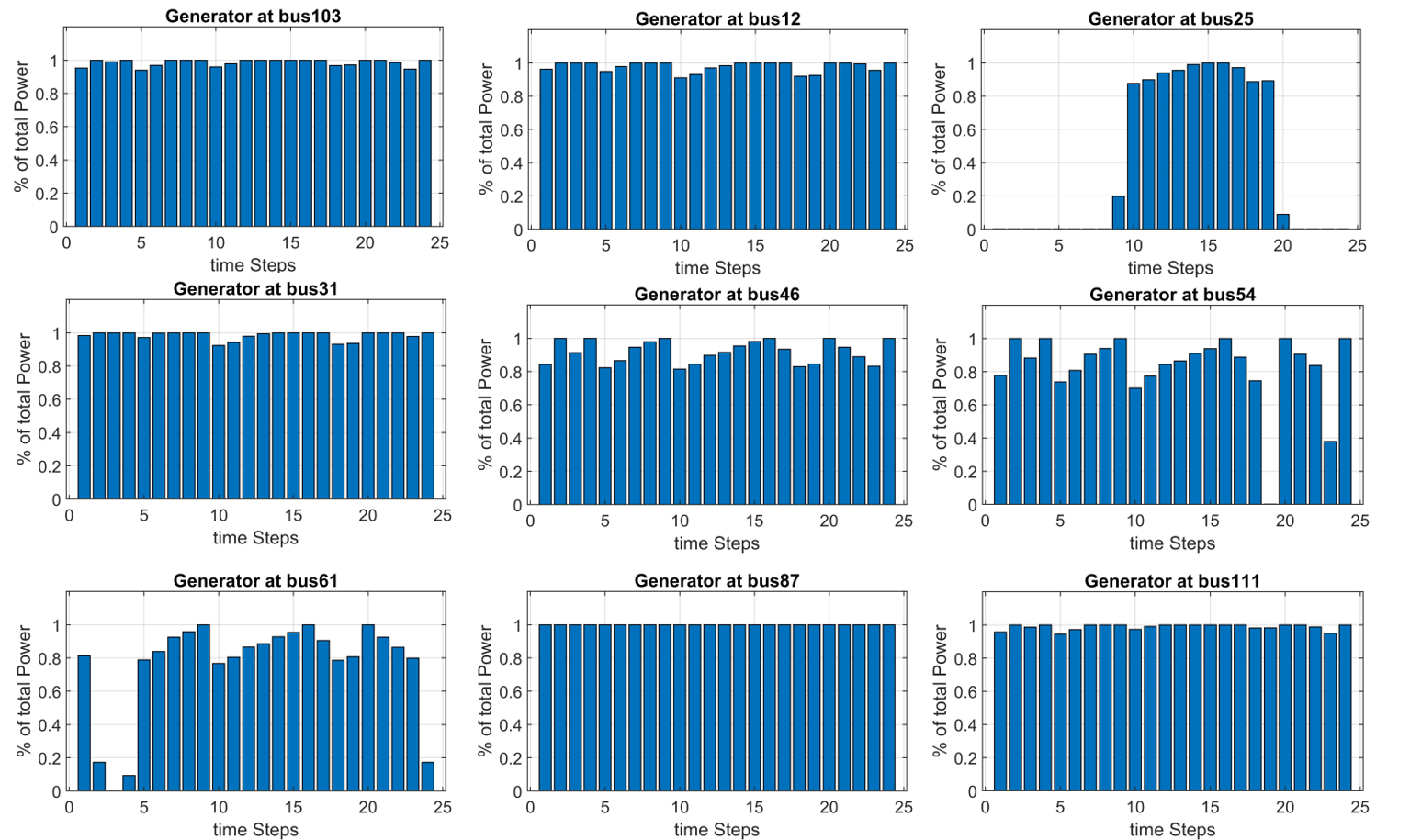}
\caption{Generator output for IEEE 118-bus system}
\label{fig:gen_118}
\end{figure}

In this study, the IEEE 118-bus system is the largest system where the system efficiency model is applied. The model was run successfully with a reasonable convergence time of 316.22s. Fig. \ref{fig:gen_118} represents the output power of the generators. For the IEEE 118-bus system, only 9 generators generated power. Only the generators with output are shown in the figure.  

\begin{figure}[t]
\centering
\includegraphics[scale=0.05]{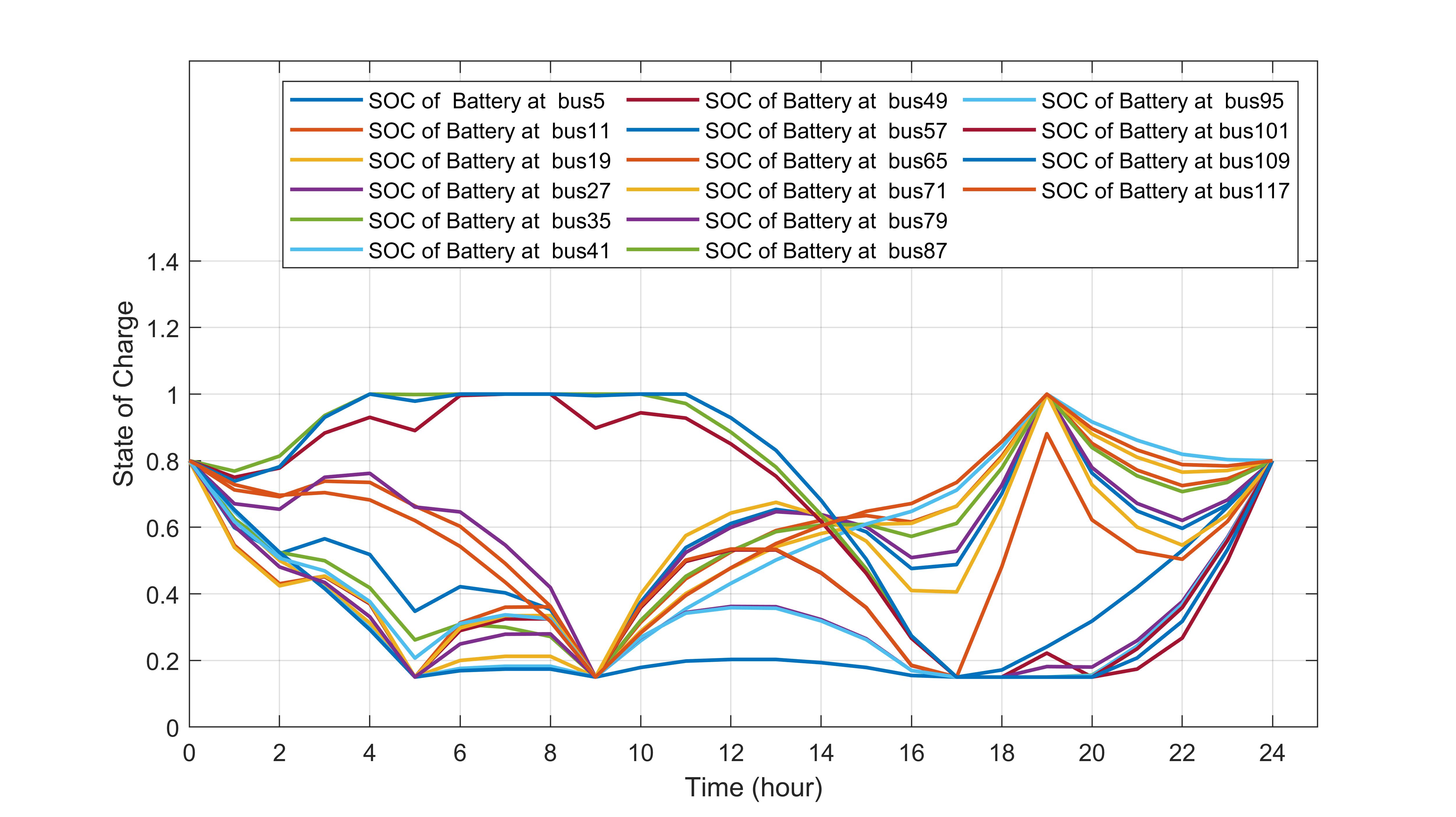}
\caption{State of charge (SOC) for IEEE 118-bus system}
\label{fig:soc_118}
\end{figure}

\begin{figure}[t]
\centering
\includegraphics[scale=0.5]{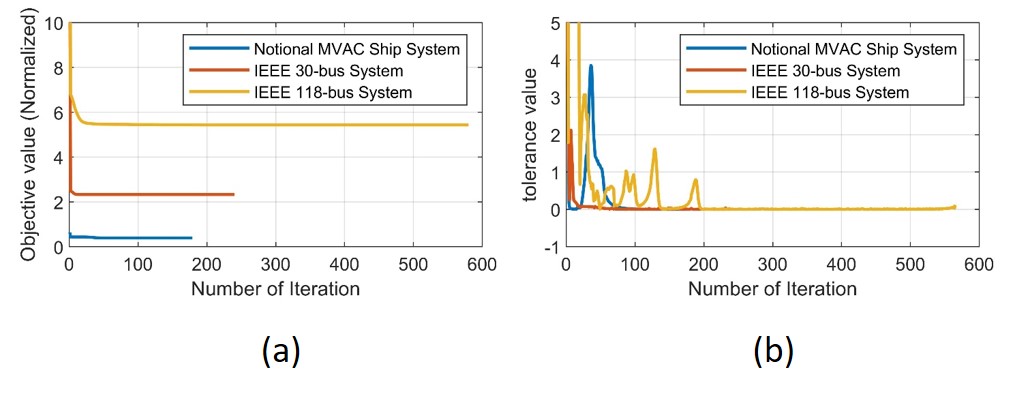}
\caption{Convergence curves of (a) Objective value (b) Tolerance value}
\label{fig:convergence}
\end{figure}

Convergence curves of both objective and tolerance values are shown in  Fig. \ref {fig:convergence}. The number of iterations required to converge depends on the number of variables in the system. Since the IEEE 118-bus system has the highest number of variables among the tested systems, it takes the most iteration to converge. 

The summary of the results for all systems is shown in TABLE \ref{tab:result_summary}. The simulations were run on Intel Core i7-10700 CPU, 2.90 GHz processor with 32.0 GB RAM. 

\begin{table}[b]
	\renewcommand{\arraystretch}{1.3}
	\caption{Summary results for notional MVAC Ship System, IEEE 30-bus, and IEEE 118-bus system}
	\label{tab:result_summary}
	\centering
	\begin{tabular}
		{
			>{\centering\arraybackslash}m{0.12\textwidth}
			>{\centering\arraybackslash}m{0.08\textwidth}
			>{\centering\arraybackslash}m{0.08\textwidth}
			>{\centering\arraybackslash}m{0.1\textwidth}
		}
		\hline
		System 
		& Number of Iteration
		& Simulation Time (s) 
		& Fuel Consumption (L) \\
		\hline
		MVAC ship system & 172 & 34.414 & $5.91\times10^4$\\
        IEEE 30-bus system & 240 & 91.728 & $3.53\times10^5$\\
        IEEE 118-bus system & 581 & 316.221 & $2.28\times10^6$\\
		\hline\hline
        \end{tabular}
\end{table}

\subsection{Performance Comparison}

In this section, a comparative analysis is conducted of the performance of the proposed fuel efficiency model in two distinct domains: 

\begin {itemize}

\item Comparison in terms of convergence time: This comparison indicates that the model proposed in this paper can be solved efficiently within a reasonable convergence time. As a result, the model can be applied to large-scale systems where the nonlinear programming model takes excessive time to converge.
\item Comparison in terms of fuel consumption: The proposed model consumes significantly less amount of fuel compared to the other models with different generator dispatches. The results indicate that the proposed model is the most efficient and optimal. \\

\end {itemize}

\subsubsection{Comparison in terms of convergence time}
The convergence time of NLP and FP models is compared in this subsection. The nonlinear optimization model from section \ref{System_Efficiency:section2}(A) is solved with MATLAB NLP (\emph{fmincon} function) for the notional MVAC ship system. The solution took an extensive time (more than 8 hours) to converge for the notional MVAC ship system with 24-time steps where the convergence time of the FP model was only 34.414s for the same system. As a result, the NLP makes the solution procedure impractical to apply to extensive systems. Moreover, the fuel consumed during the operation was $6.53\times10^4$L, higher than the fuel consumed with the proposed FP model. The FP model is clearly more advantageous than the nonlinear programming optimization, even with a higher number of variables. 
The performance comparison between FP and NLP is shown in TABLE \ref{tab:FP_vs_NLP}.

\begin{table} [b]
    \centering
    \renewcommand{\arraystretch}{1.8}
    \caption{Performance comparison between FP and NLP for notional MVAC ship system}
    \label{tab:FP_vs_NLP}
    \begin{tabular}
    {
			>{\centering\arraybackslash}m{0.08\textwidth}
			>{\centering\arraybackslash}m{0.07\textwidth}
			>{\centering\arraybackslash}m{0.06\textwidth}
			>{\centering\arraybackslash}m{0.07\textwidth}
            >{\centering\arraybackslash}m{0.07\textwidth}
		}
        \hline \hline

        System & Optimization Model & Total Variables & Simulation Time (s) & Fuel Consumption (L) \\
        \hline
        \multirow{2}{5em}{Notional MVAC Ship System} & FP & 2592 & 34.414 & $5.91\times10^4$ \\

        & NLP & 1488 & 30233.884 & $6.53\times 10^4$ \\
        \hline
        \multirow{2}{*}{IEEE 30-bus} & FP & 6960 & 91.728 & $3.53\times10^5$ \\
 
        & NLP & 4128 & did not converge & N\slash A \\
        \hline
        \multirow{2}{*}{IEEE 118-bus} & FP & 28128 & 316.221 & $2.28\times10^6$ \\
 
        & NLP & 17088 & did not converge & N\slash A \\
        \hline \hline
    \end{tabular}
\end{table}

\subsubsection{Comparison in terms of fuel consumption}
The proposed system efficiency model is compared with three other models (as listed in TABLE \ref{tab:models}) to demonstrate the fuel efficiency. The difference between the models is in the generator dispatch allowed during the simulation. The number and type of generators allowed for each model during the operation are indicated by the '\emph{Dispatch}' column. 

\begin{figure}[t]
\centering
\includegraphics[scale=0.1]{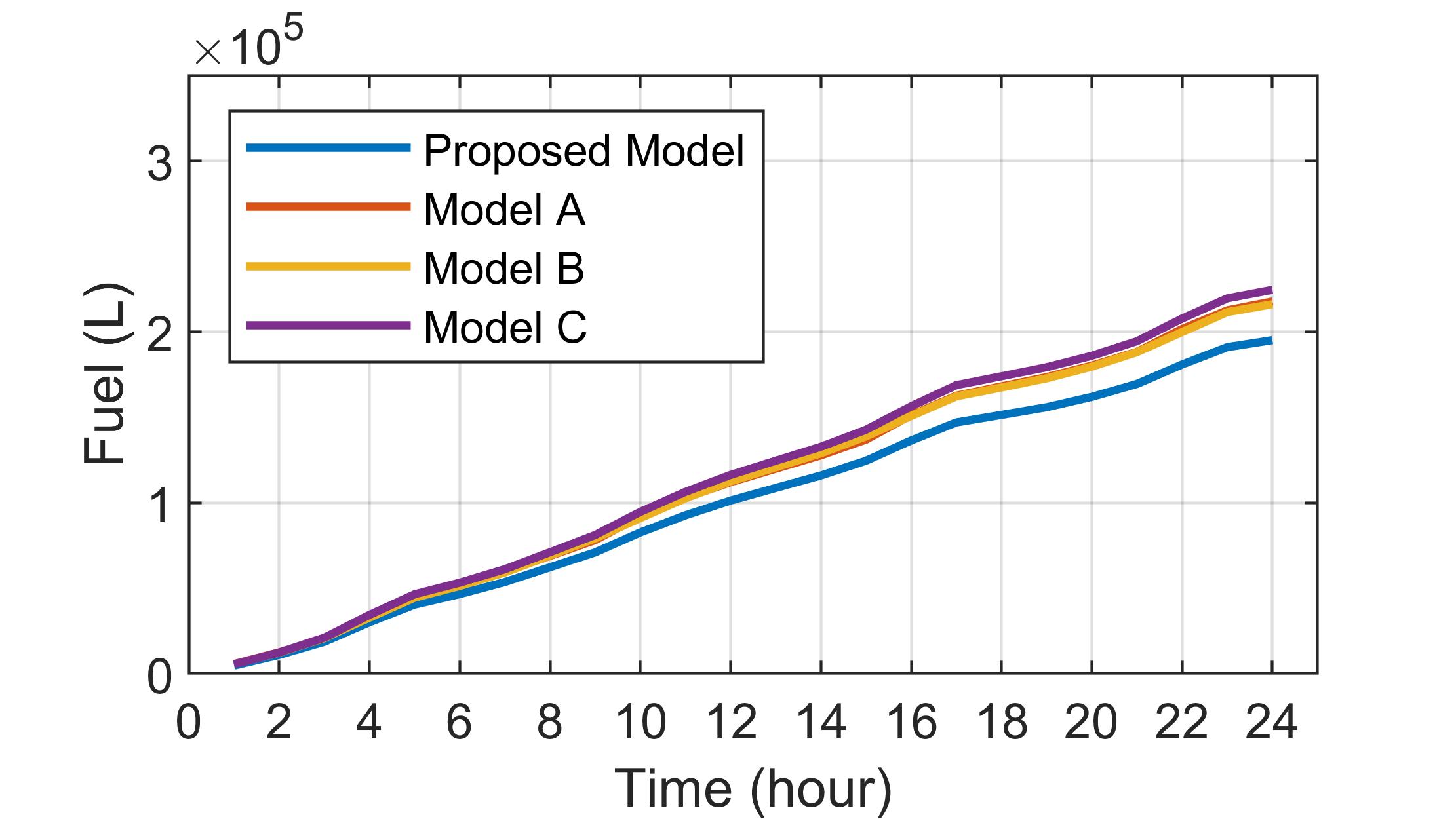}
\caption{Fuel consumption comparison between different models of notional MVAC system}
\label{fig:fuel_consumption}
\end{figure}

The comparison results are shown in Fig. \ref{fig:fuel_consumption}. It can be observed that the proposed model has the lowest fuel consumption among all models. Initially, all models consume almost similar amounts of fuel. However, the difference in fuel consumption increases with the number of time steps. This observation highlights the superior efficiency of the proposed model in comparison to the models examined in this section.

\begin{table}
	\renewcommand{\arraystretch}{1.2}
	\caption{Summary of generator dispatch of different models of notional MVAC ship system}
	\label{tab:models}
	\centering
	\begin{tabular}
		{
			>{\centering\arraybackslash}m{0.18\textwidth}
			>{\centering\arraybackslash}m{0.25\textwidth}
		}
		\hline
		Models 
		& Dispatch \\
		\hline
		Proposed Model & 2 MTGs and 2 ATGs\\
        Model A & 2 MTGs\\
        Model B & 1 MTG and 1 ATG \\
        Model C & Equal power-sharing \\
		\hline\hline
        \end{tabular}
\end{table}

\section{Conclusion} \label{System_Efficiency:section5}
This study has addressed the challenge of the fuel consumption minimization problem to enhance the system efficiency and reduce the operating cost of the power generation units. The traditional approaches typically focus on maximizing the efficiency function or minimizing the generator cost function to achieve optimal fuel consumption for the system. However, these approaches do not account for the fuel consumption rate directly and are impractical to implement in real-world systems where optimizing fuel use is the objective. In addition, existing studies that have used the fuel consumption curve to formulate the optimization problem have numerous limitations, including incompatibility to apply to large AC systems.

As a result, it is crucial to incorporate a function that directly represents fuel consumption to enhance the system's fuel efficiency. This study introduced a novel objective function based on a sum-of-ratios approach, providing a straightforward representation of the fuel consumption rate. The sum-of-ratios problem was effectively solved by leveraging a fractional programming (FP) reformulation technique, resulting in successful fuel consumption minimization. Moreover, the low convergence time of the solution makes the model suitable for large-scale systems. While the model stands out in its uniqueness and effectiveness compared to other approaches, future research will concentrate on implementing a distributed algorithm to enhance scalability for larger and more complex systems. 

\section{Acknowledgement} \label{System_Efficiency:section6}
The information, data, or work presented herein was partly funded by the U.S. Office of Naval Research under the award numbers N000142212239 and N000142112124.

\bibliographystyle{IEEEtran}
\bibliography{references}

\end{document}